\documentclass[oneside, reqno] {amsart}

\usepackage{xypic}
\input xy
\xyoption{all}
\usepackage{epsfig}
\usepackage{amsthm}
\usepackage{amssymb}
\usepackage{amsmath}
\usepackage{amscd}
\usepackage{color}
\usepackage[T1]{fontenc}
\usepackage{upgreek}
\usepackage[font=scriptsize]{caption}
\usepackage{wrapfig}
\usepackage{graphicx}
\usepackage{subfigure}

\newcommand{\bg}{\begin{equation}}
\newcommand{\ed}{\end{equation}}
\newcommand{\bga}{\begin{eqnarray}}
\newcommand{\eda}{\end{eqnarray}}

\def\cbdu{\par{\raggedleft$\Box$\par}}

\newtheorem {Theorem}  {Theorem}

\numberwithin{Theorem}{section}

\newtheorem {Lemma}[Theorem]  {Lemma}

\theoremstyle{definition}

\theoremstyle{remark}
\newtheorem{Remark}[Theorem]{\bf Remark}

\expandafter\chardef\csname pre amssym.def
at\endcsname=\the\catcode`\@ \catcode`\@=11
\def\undefine#1{\let#1\undefined}
\def\newsymbol#1#2#3#4#5{\let\next@\relax
 \ifnum#2=\@ne\let\next@\msafam@\else
 \ifnum#2=\tw@\let\next@\msbfam@\fi\fi
 \mathchardef#1="#3\next@#4#5}
\def\mathhexbox@#1#2#3{\relax
 \ifmmode\mathpalette{}{\m@th\mathchar"#1#2#3}%
 \else\leavevmode\hbox{$\m@th\mathchar"#1#2#3$}\fi}
\def\hexnumber@#1{\ifcase#1 0\or 1\or 2\or 3\or 4\or 5\or 6\or 7\or 8\or
 9\or A\or B\or C\or D\or E\or F\fi}

\font\teneufm=eufm10 \font\seveneufm=eufm7 \font\fiveeufm=eufm5
\newfam\eufmfam
\textfont\eufmfam=\teneufm \scriptfont\eufmfam=\seveneufm
\scriptscriptfont\eufmfam=\fiveeufm

\catcode`\@=\csname pre amssym.def at\endcsname

\newcounter{remark}
\def  \12  {{\frac{1}{2}}}

\def\build#1_#2^#3{\mathrel{\mathop{\kern 0pt#1}\limits_{#2}^{#3}}}

\numberwithin{equation}{section}

\begin{document}
%\currannalsline{0}{2006}

%\title[Local existence for Hall-MHD]{Local existence for the non-resistive Hall-magneto-hydrodynamics system in $\R^n$}
\title[Dissipation wavenumber and regularity]{Dissipation wavenumber and regularity for electron magnetohydrodynamics}

%\author{hello}

\author [Mimi Dai]{Mimi Dai}

\address{Department of Mathematics, Statistics and Computer Science, University of Illinois at Chicago, Chicago, IL 60607, USA}
\email{mdai@uic.edu}

\author [Chao Wu]{Chao Wu}

\address{Department of Mathematics, Statistics and Computer Science, University of Illinois at Chicago, Chicago, IL 60607, USA}
\email{cwu206@uic.edu}

\thanks{The authors are partially supported by the NSF grant DMS--2009422. M.Dai also acknowledges the support of the AMS Centennial Fellowship.}

\begin{abstract}
We consider the electron magnetohydrodynamics (MHD) with static background ion flow. A special situation of $B(x,y,t)=\nabla\times (a\vec e_z)+b \vec e_z$ with scalar-valued functions $a(x,y,t)$ and $b(x,y,t)$ was studied numerically in the physics paper \cite{WH}. The authors concluded from numerical simulations that there is no evidence of dissipation cutoff for the electron MHD. In this paper we show the existence of determining wavenumber for the electron MHD, and establish a regularity condition only on the low modes of the solution.
Our results suggest that the conclusion of the physics paper on the dissipation cutoff for the electron MHD is debatable.

\bigskip


\hspace{0.02cm}CLASSIFICATION CODE: 35Q35, 65M70, 76D03, 76W05, 76B03.
\end{abstract}

\maketitle

\section{Introduction}
As an important mathematical model to describe the fast magnetic reconnection phenomena in plasma physics, the incompressible magnetohydrodynamics (MHD) system with Hall effect 
\begin{equation}\label{mhd}
\begin{split}
u_t+(u\cdot\nabla) u-(B\cdot\nabla) B+\nabla P=&\ \nu\Delta u, \\
B_t+(u\cdot\nabla) B-(B\cdot\nabla) u 
+\nabla\times ((\nabla\times B)\times B)=&\ \mu\Delta B, \\
\nabla\cdot u=0, \ \ \nabla\cdot B=&\ 0,
\end{split}
\end{equation}
derived in \cite{ADFL} differs from the classical MHD system by including the additional Hall term $\nabla\times ((\nabla\times B)\times B)$. As customary, the unknowns of (\ref{mhd}) are the velocity field $u$, the magnetic field $B$ and the fluid pressure $P$. The parameters $\nu$ and $\mu$ stand for the kinetic viscosity and magnetic resistivity respectively. The Hall term is notoriously more challenging than other nonlinear terms in the system. With the aim to have a thorough understanding of the nonlinear structure of the Hall term we consider the electron MHD 
\begin{equation}\label{emhd}
\begin{split}
B_t+ \nabla\times ((\nabla\times B)\times B)=&\ \mu\Delta B,\\
\nabla\cdot B=&\ 0
\end{split}
\end{equation}
by taking $u\equiv 0$ in (\ref{mhd}). The electron MHD recounts the magnetohydrodynamics with static background ion flow. It has been studied by both mathematicians and physicists, for instance see \cite{JO, WH}.
We observe the first equation of (\ref{emhd}) is quasi-linear, which sets immediate obstacles for mathematical analysis of the model. 

In \cite{WH} the authors studied the electron MHD under a two dimensional setting, that is
\begin{equation}\label{two-half}
B(x,y,t)=\nabla\times (a\vec e_z)+b\vec e_z
\end{equation}
with scalar-valued functions 
\[a=a(x,y,t), \ \ b=b(x,y,t),\]
and $\vec e_z=(0, 0,1)$. Note that $\nabla\cdot B=0$ and $B=(a_y, -a_x, b)$. Moreover, it follows from (\ref{emhd}) that $a$ and $b$ satisfy the system
\begin{equation}\label{emhd1}
\begin{split}
a_t+ (a_yb_x-a_xb_y)=&\ \mu\Delta a, \\
b_t-(a_y\Delta a_x-a_x\Delta a_y)=&\ \mu\Delta b.
\end{split}
\end{equation}
We observe that the equation of $a$ is linear in $a$, so is the equation of $b$. However it is not necessarily to indicate that system (\ref{emhd1}) is easy to be analyzed. In fact, one can see system (\ref{emhd1}) is supercritical via the following
argument. Formally for smooth functions $a$ and $b$, the basic energy law for (\ref{emhd1}) is given by
\begin{equation}\label{energy}
\frac12\frac{d}{dt}\int \left(a_x^2+a_y^2+b^2\right)\, d\vec x+\mu \int \left(a_{xx}^2+2a_{xy}^2+a_{yy}^2+b_x^2+b_y^2\right)\, d\vec x=0
\end{equation}
with $\vec x=(x,y,0)$. It gives the a priori energy estimates 
\begin{equation}\label{basic-space}
\begin{split}
& a\in L^\infty(0,T; H^1(\mathbb T^2))\cap L^2(0,T; H^2(\mathbb T^2)), \\ 
& b\in L^\infty(0,T; L^2(\mathbb T^2))\cap L^2(0,T; H^1(\mathbb T^2)).
\end{split}
\end{equation}
On the other hand, system (\ref{emhd1}) has the natural scaling
\[a_\lambda=\lambda^{-1} a(\lambda \vec x, \lambda^2t), \ \ b_\lambda= b(\lambda \vec x, \lambda^2t).\]
In two spatial dimension, some critical spaces for $a$ are 
\[\dot H^2\subset \dot B^{1+\frac2p}_{p, \infty}\subset \dot B^{1}_{\infty, \infty}, \ \ 2\leq p<\infty, \]
and some critical spaces for $b$ are
\[\dot H^1\subset \dot B^{\frac2p}_{p, \infty}\subset \dot B^{0}_{\infty, \infty}, \ \ 2\leq p<\infty.\]
Thus in view of (\ref{basic-space}), system (\ref{emhd1}) is energy supercritical. In general, global existence of regular solution of (\ref{emhd1}) is not expected. 

Numerical simulations for system (\ref{emhd1}) were performed in \cite{WH}. Two main discoveries of \cite{WH} are: (i) no evidence for dissipation cutoff is found; (ii) stationary structures are found to develop in time. 
In this paper we disprove (i) by showing the existence of determining wavenumber and establishing a low modes regularity criterion. In order to verify (ii), investigation on steady states and stability will be followed in future work. 

We also make the following observations. 
If $a$ is a radial function, that is, $a(x,y,t)=f(x^2+y^2, t)$ for some function $f$, we have
\[a_y\Delta a_x-a_x\Delta a_y=0.\] 
Hence in this case system (\ref{emhd1}) becomes
\begin{equation}\label{emhd2}
\begin{split}
a_t+ (a_yb_x-a_xb_y)=&\ \mu\Delta a, \\
b_t=&\ \mu\Delta b
\end{split}
\end{equation}
where $b$ satisfies a heat equation. If both $a$ and $b$ are radial functions, we further have $a_yb_x-a_xb_y=0$. Therefore for radial functions $a$ and $b$, the magnetic field $B$ given by (\ref{two-half}) is a steady state of the ideal electron MHD since 
\[\nabla\times ((\nabla\times B)\times B)=0.\]

\medskip

\subsection{Main results}

We first consider the general form (\ref{emhd}) with a given external forcing $f$ on $\mathbb T^n$ with $n=2,3$,
\begin{equation}\label{emhd-f}
\begin{split}
B_t+ \nabla\times ((\nabla\times B)\times B)=&\ \mu\Delta B+f,\\
\nabla\cdot B=&\ 0.
\end{split}
\end{equation}
When $n=2$, it is understood that
\[B=B(x,y,t)=\left(B_1(x,y,t), B_2(x,y,t), B_3(x,y,t) \right).\]
We show the existence of determining wavenumber for (\ref{emhd-f}) following the wavenumber splitting framework developed in \cite{CD-nse-modes, CDK}. 

According to the scaling of (\ref{emhd-f}), we define the wavenumber 
\begin{equation}\label{wave}
\lambda_{Q(B)}(t)=\min \{\lambda_q: \lambda_p^{\frac{n}{r}}\|B_p\|_{L^r}<c_r \mu, \forall \ p>q; \ \mbox{and} \ \ \|B_{\leq q}\|_{L^\infty}<c_r\mu, q\in \mathbb N\}
\end{equation}
for some parameter $r\in (n, 2n)$ and a small constant $c_r>0$. Note both quantities $\lambda_p^{\frac{n}{r}}\|B_p\|_{L^r}$ and $\|B_{\leq q}\|_{L^\infty}$ are scaling invariant. To simplify notations, $c_r$ will be chosen small enough such that it absorbs all the other constants appeared in estimates in Section \ref{sec-det}. We prove

\begin{Theorem}\label{thm-det}
Let $B^{(1)}(t)$ and $B^{(2)}(t)$ be two solutions of (\ref{emhd}). Let $\lambda_{Q(B^{(1)})}(t)$ and $\lambda_{Q(B^{(2)})}(t)$ be the wavenumber defined for $B^{(1)}(t)$ and $B^{(2)}(t)$ respectively as in (\ref{wave}). Denote 
\[\lambda_{Q}(t)=\max\{\lambda_{Q(B^{(1)})}(t), \lambda_{Q(B^{(2)})}(t)\}.\]
Assume \[ \left(B^{(1)}-B^{(2)}\right)|_{\leq \lambda_{Q}}=0.\]
Then we have 
\[\lim_{t\to\infty} \|B^{(1)}-B^{(2)}\|_{H^s}=0, \ \ \ \mbox{for} \ \ -\frac{n}{r}<s<\frac{n}{r}-1.\]
\end{Theorem}
Since $r\in (n,2n)$, $s$ can be as close as to $0$; and hence the convergence occurs in space close to $L^2(\mathbb T^n)$.

For the particular form (\ref{emhd1}) numerically studied in \cite{WH}, we further establish the following regularity condition only on the low modes of the solution. 
Define 
\begin{equation}\label{wave-a}
\lambda_{Q(a)}(t)=\min\left\{ \lambda_q: \lambda_p^{\frac{2}{r}}\|\nabla\times (a_p\vec e_z)\|_{L^r}<c_r \mu \ \ \forall p>q\right\},
\end{equation}
\begin{equation}\label{wave-b}
\lambda_{Q(b)}(t)=\min\left\{ \lambda_q: \lambda_p^{\frac{2}{r}}\|b_p\|_{L^r}<c_r \mu \ \ \forall p>q\right\},
\end{equation}
for some $r\in[2,\infty)$ and a small enough constant $c_r$. 
Denote 
\[f_1(t)= \|\nabla b_{\leq Q(b)}(t)\|_{B^1_{\infty,\infty}}, \ \ f_2(t)= \|\nabla \nabla a_{\leq Q(a)}(t)\|_{B^1_{\infty,\infty}}.\]

\begin{Theorem}\label{thm-reg-low}
Assume $(a(t),b(t))$ is a regular solution of (\ref{emhd1}) on $[0,T)$ and $f_1, f_2\in L^1(0,T)$. Then $(a(t),b(t))$ is regular on $[0,T]$.  
\end{Theorem}
The result of Theorem \ref{thm-reg-low} is consistent with the regularity criterion obtained for the 3D electron MHD in \cite{Dai-hmhd-reg}. A modification of the proof for the 3D case from \cite{Dai-hmhd-reg} would work for the 2D case considered here. However we provide a self-contained proof for Theorem \ref{thm-reg-low} in Section \ref{sec-reg}, by exploring cancellations involving the interactions of the horizontal component $\nabla\times (a\vec e_z)$ and the vertical component $b\vec e_z$. Such cancellations are believed to be beneficial for our future investigation of (\ref{emhd1}). 

We make the observation that both results in Theorem \ref{thm-det} and Theorem \ref{thm-reg-low} suggest the existence of dissipation cutoff for the electron MHD. Therefore the discovery of no evidence for dissipation cutoff based on numerical simulation in \cite{WH} has no theoretical ground.

%{\color{blue} Question: How to make sure $a(t)$ remains a radial function? }

\medskip

%In the general case when $a$ is not necessarily radial, we have the conditional regularity result for (\ref{emhd1}) with conditions imposed only on $b$. 
For (\ref{emhd1}) we also provide a regularity condition of the Ladyzhenskaya-Prodi-Serrin type only on the vertical component $b\vec e_z$.
\begin{Theorem}\label{thm-criterion}
Let $(a(t),b(t))$ be a regular solution to (\ref{emhd1}) on $[0,T)$ for any $T>0$. If 
\begin{equation}\label{criterion}
%\int_0^T \|\nabla b\|_{L^\infty}^2\, dt<\infty,
b\in L^s(0,T; L^r(\mathbb T^2)) \ \ \mbox{with} \ \ \frac2s+\frac2r\leq 1 \ \ \forall \ \ r\in(2, \infty],
\end{equation}
the solution $(a(t),b(t))$ is regular on $[0,T]$.
\end{Theorem}

The regularity criterion (\ref{criterion}) for the general 2D electron MHD was obtained recently in \cite{RY} where the authors discovered some hidden cancellations. In the case of (\ref{emhd1}), the cancellations are more transparent. Again, such cancellations may play crucial roles in future study of (\ref{emhd1}). Therefore we present our own proof of \ref{thm-criterion} in Section \ref{sec-criterion}. 

In the end we make the observation that the solution of (\ref{emhd1}) with radial horizontal component is regular. 

\begin{Theorem}\label{thm-global}
Let $a_0\in H^2(\mathbb T^2)$ and $b_0\in H^1(\mathbb T^2)$. Assume $a_0$ is a radial function, i.e. $a_0=a_0(r)$ with $r=(x^2+y^2)^{\frac12}$. There exists a global regular solution $(a(t), b(t))$ to system (\ref{emhd2}) with the initial data $(a_0, b_0)$. 
\end{Theorem}

Note that the proof of Theorem \ref{thm-global} follows from Theorem \ref{thm-criterion} immediately. Indeed, $b$ satisfies the heat equation in system (\ref{emhd2}) and hence the condition (\ref{criterion}) with $s=2$ and $r=\infty$ is satisfied.

\bigskip

\section{Preliminaries}

\subsection{Notations}
A general constant which is not important to be tracked is denoted by $C$, which may be different from line to line.  The notation $\lesssim$ is used to replace $\leq$ up to a multiplication of constant.

\medskip

\subsection{Littlewood-Paley theory} To apply the frequency localization techniques,  the basics of Littlewood-Paley theory on $\mathbb T^n=[-L, L]^n$ are recalled briefly.  Let $\lambda_q=2^q/L$ for any integer $q$. Define the radial function $\chi\in C_0^\infty(\mathbb R^n)$ as
\begin{equation}\notag
\chi(\xi)=
\begin{cases}
1, \ \ \ \mbox{for} \ \ |\xi|\leq \frac34\\
0, \ \ \ \mbox{for} \ \ |\xi|\geq 1.
\end{cases}
\end{equation}
Denote $\varphi(\xi)=\chi(\frac{\xi}{2})-\chi(\xi)$ and 
\begin{equation}\notag
\varphi_q(\xi)=
\begin{cases}
\varphi(\lambda_q^{-1}\xi), \ \ \ \mbox{for} \ \ q\geq0\\
\chi(\xi), \ \ \ \ \ \ \ \  \mbox{for} \ \ q= -1. 
\end{cases}
\end{equation}
The $q$-th Littlewood-Paley projection of a tempered distribution vector field $u$ on $\mathbb T^n$ is defined by 
\[u_q(x)=\Delta_q u(x):= \sum_{k\in \mathbb Z^n} \widehat u(k)\varphi_q(k) e^{i\frac{2\pi}{L} k\cdot x}\]
where $\widehat u(k)$ is the Fourier coefficient of $u$.  The unit decomposition 
\[u=\sum_{q=-1}^\infty u_q\]
holds in the distributional sense. Denote
\[u_{\leq Q}=\sum_{q\geq-1}^Q u_q, \ \ \ \ \ \ \widetilde {\Delta_q} u=\widetilde u_q=u_{q-1}+u_q+u_{q+1}.\]
It is convenient to use the equivalent norm of $u$ in the space $H^s(\mathbb T^n)$ 
\[\|u\|_{H^s(\mathbb T^n)}=\left(\sum_{q=-1}^\infty \lambda_q^{2s} \|u_q\|_{L^2}^2\right)^{\frac12}.\]

We will use the following Bony's paraproduct formula 
\begin{equation}\notag
\begin{split}
\Delta_q(u\cdot v)=&\sum_{|p-q|\leq 2} \Delta_q(u_{\leq p-2}\cdot \nabla v_p)+\sum_{|p-q|\leq 2} \Delta_q(u_{p}\cdot \nabla v_{\leq p-2})\\
&+\sum_{p\geq q-2} \Delta_q(\widetilde u_p\cdot \nabla v_p),
\end{split}
\end{equation}
which applies to scalar-valued functions as well. 
We recall some commutator estimates from \cite{Dai-hmhd-reg}. 
Denote 
\[[\Delta_q, u_{\leq p-2}\cdot\nabla ]v_p=\Delta_q(u_{\leq p-2}\cdot \nabla v_p)-u_{\leq p-2}\cdot \nabla \Delta_q v_p,\]
\[[\Delta_q, u\times \nabla\times ]v=\Delta_q(u\times(\nabla\times v)) -u\times(\nabla\times v_q).\]

\begin{Lemma}\label{le-comm1}
For any $1<r_1<\infty$ and $2\leq r_2\leq \infty$ with $\frac1{r_1}=\frac1{r_2}+\frac1{r_3}$, we have
\[ \|[\Delta_q, u_{\leq p-2}\cdot\nabla ]v_p\|_{L^{r_1}} \lesssim \|\nabla u_{\leq p-2}\|_{L^{r_2}}\| v_{p}\|_{L^{r_3}}, \]
\[ \|[\Delta_q, u\times \nabla\times ]v\|_{L^{r_1}} \lesssim \|\nabla u\|_{L^{r_2}}\| v\|_{L^{r_3}}, \]
\end{Lemma}

\begin{Lemma}\label{le-comm2}
For any $1<r_1, r_2<\infty$ with $\frac1{r_1}+\frac1{r_2}=1$, the estimate
\[ \int_{\mathbb T^n} [\Delta_q, u\times \nabla\times ]v\cdot \nabla\times w\, d\vec x\lesssim \|\nabla\nabla u\|_{L^\infty} \| v\|_{L^{r_1}}\| w\|_{L^{r_2}}.\]
holds. 
\end{Lemma}
The above commutators and commutator estimates also apply to scalar-valued functions with the form of the differential operators modified appropriately.

\bigskip

\section{Determining wavenumber}
\label{sec-det}

We prove Theorem \ref{thm-det} in this section. 
Denote $h(t)=B^{(1)}(t)-B^{(2)}(t)$, which satisfies the equation
\begin{equation}\label{eq-H}
h_t+\nabla\times \left((\nabla\times B^{(1)})\times h\right)+\nabla\times \left((\nabla\times h)\times B^{(2)}\right)=\mu\Delta h
\end{equation}
with $\nabla \cdot h=0$. It is clear from the assumption of the theorem that $h|_{\leq Q}\equiv 0$. On the other hand, thanks to the definition of $\lambda_Q$ we have 
\begin{equation}\label{high-modes}
\begin{split}
\lambda_q^{\frac{n}{r}}\|B^{(1), (2)}_q\|_{L^r}<&\ c_r \mu, \ \ \forall \ \ q>Q, \\
\|B^{(1), (2)}_{\leq Q}\|_{L^\infty}<&\ c_r \mu.
\end{split}
\end{equation}

We estimate the $H^{s}$ norm of $h(t)$ for appropriate $s$ in the following. It follows from (\ref{eq-H})
\begin{equation}\label{est-energy1}
\begin{split}
&\frac12\frac{d}{dt} \sum_{q\geq -1} \lambda_q^{2s}\|h_q\|_{L^2(\mathbb T^n)}^2+\mu \sum_{q\geq -1} \lambda_q^{2s+2}\|h_q\|_{L^2(\mathbb T^n)}^2\\
=&-\sum_{q\geq -1} \lambda_q^{2s} \int_{\mathbb T^n} \Delta_q\left((\nabla\times B^{(1)})\times h\right)\cdot \nabla\times h_q\, d\vec x\\
&-\sum_{q\geq -1} \lambda_q^{2s} \int_{\mathbb T^n} \Delta_q\left((\nabla\times h)\times B^{(2)}\right)\cdot \nabla\times h_q\, d\vec x\\
=&: -I-J.
\end{split}
\end{equation}

The estimate of $I$ starts with the decomposition of Bony's paraproduct 
\begin{equation}\notag
\begin{split}
I=&\sum_{q\geq -1} \sum_{|p-q|\leq 2} \lambda_q^{2s} \int_{\mathbb T^n} \Delta_q\left((\nabla\times B^{(1)}_{\leq p-2})\times h_p\right)\cdot \nabla\times h_q\, d\vec x\\
&+\sum_{q\geq -1}\sum_{|p-q|\leq 2} \lambda_q^{2s} \int_{\mathbb T^n} \Delta_q\left((\nabla\times B^{(1)}_p)\times h_{\leq p-2}\right)\cdot \nabla\times h_q\, d\vec x\\
&+\sum_{q\geq -1} \sum_{p\geq q-2}\lambda_q^{2s} \int_{\mathbb T^n} \Delta_q\left((\nabla\times \widetilde B^{(1)}_p)\times h_p\right)\cdot \nabla\times h_q\, d\vec x\\
=&: I_1+I_2+I_3. 
\end{split}
\end{equation}
Since $h_{\leq Q}\equiv 0$, we split $I_1$ as 
\begin{equation}\notag
\begin{split}
I_1=& \sum_{q>Q}\sum_{\substack {|p-q|\leq 2\\ p>Q+2}}\lambda_q^{2s} \int_{\mathbb T^n} \Delta_q\left((\nabla\times B^{(1)}_{(Q, p-2]}\times h_p\right)\cdot \nabla\times h_q\, d\vec x\\
& +\sum_{q>Q}\sum_{\substack {|p-q|\leq 2\\ p>Q+2}}\lambda_q^{2s} \int_{\mathbb T^n} \Delta_q\left((\nabla\times B^{(1)}_{\leq Q}\times h_p\right)\cdot \nabla\times h_q\, d\vec x\\
=&: I_{11}+I_{12}.
\end{split}
\end{equation}
Applying H\"older's inequality, Bernstein's inequality and (\ref{high-modes}) gives
\begin{equation}\notag
\begin{split}
|I_{11}|\lesssim& \sum_{q>Q} \sum_{\substack {|p-q|\leq 2\\ p>Q+2}}\lambda_q^{2s}\|h_p\|_{L^2} \|\nabla\times h_q\|_{L^{\frac{2r}{r-2}}}
\sum_{Q<p'\leq p-2} \|\nabla\times B^{(1)}_{p'}\|_{L^r}\\
\leq&\ c_0\mu \sum_{q>Q} \sum_{\substack {|p-q|\leq 2\\ p>Q+2}}\lambda_q^{2s+1+\frac{n}{r}}\|h_p\|_{L^2} \|h_q\|_{L^2}
\sum_{Q<p'\leq p-2}\lambda_{p'}^{1-\frac{n}{r}}\\
\leq &\ c_0\mu \sum_{q>Q}\lambda_q^{2s+1+\frac{n}{r}} \|h_q\|_{L^2}^2\sum_{Q<p'\leq q}\lambda_{p'}^{1-\frac{n}{r}}\\
\leq &\ c_0\mu \sum_{q>Q}\lambda_q^{2s+2} \|h_q\|_{L^2}^2\sum_{Q<p'\leq q}\lambda_{p'-q}^{1-\frac{n}{r}}\\
\leq &\ c_0\mu \sum_{q>Q}\lambda_q^{2s+2} \|h_q\|_{L^2}^2
\end{split}
\end{equation}
where we used $1-\frac{n}{r}>0$ in the last step. Similarly we have
\begin{equation}\notag
\begin{split}
|I_{12}|\lesssim& \sum_{q>Q} \sum_{\substack {|p-q|\leq 2\\ p>Q+2}}\lambda_q^{2s} \|\nabla\times B^{(1)}_{\leq Q}\|_{L^\infty} \|h_p\|_{L^2}\|\nabla\times h_q\|_{L^2}\\
\leq&\ c_0\mu \sum_{q>Q} \sum_{\substack {|p-q|\leq 2\\ p>Q+2}}\lambda_q^{2s} \lambda_Q \|h_p\|_{L^2}\|\nabla\times h_q\|_{L^2}\\
\leq &\ c_0\mu \sum_{q>Q}\lambda_q^{2s+2} \|h_q\|_{L^2}^2.
\end{split}
\end{equation}
The term $I_2$ is handled analogously,
\begin{equation}\notag
\begin{split}
|I_2|\leq & \sum_{q>Q} \sum_{\substack {|p-q|\leq 2\\ p>Q+2}}\lambda_q^{2s} \int_{\mathbb T^n} \left|\Delta_q\left((\nabla\times B^{(1)}_p)\times h_{(Q, p-2]}\right)\cdot \nabla\times h_q\right|\, d\vec x\\
\leq & \sum_{q>Q} \sum_{\substack {|p-q|\leq 2\\ p>Q+2}}\lambda_q^{2s}\|\nabla\times B^{(1)}_p\|_{L^r}\|h_{(Q,p-2]}\|_{L^{\frac{2r}{r-2}}}\|\nabla\times h_q\|_{L^2}\\
\leq &\ c_0\mu  \sum_{q>Q} \sum_{\substack {|p-q|\leq 2\\ p>Q+2}}\lambda_q^{2s+1}\lambda_p^{1-\frac{n}{r}}\|h_q\|_{L^2}\|h_{(Q,p-2]}\|_{L^{\frac{2r}{r-2}}}\\
\leq &\ c_0\mu  \sum_{q>Q}\lambda_q^{2s+2-\frac{n}{r}}\|h_q\|_{L^2}\sum_{Q<p'\leq q}\lambda_{p'}^{\frac{n}{r}}\|h_{p'}\|_{L^2}\\
\leq &\ c_0\mu  \sum_{q>Q}\lambda_q^{s+1}\|h_q\|_{L^2}\sum_{Q<p'\leq q}\lambda_{p'}^{s+1}\|h_{p'}\|_{L^2} \lambda_{q-p'}^{s+1-\frac{n}{r}},
\end{split}
\end{equation}
followed by using Jensen's inequality for $s+1-\frac{n}{r}<0$
\begin{equation}\notag
\begin{split}
|I_2| \leq& c_0\mu  \sum_{q>Q}\lambda_q^{2s+2}\|h_q\|_{L^2}^2+c_0\mu \sum_{q>Q} \sum_{Q<p'\leq q}\lambda_{p'}^{2s+2}\|h_{p'}\|_{L^2}^2 \lambda_{q-p'}^{s+1-\frac{n}{r}}\\
\leq& c_0\mu  \sum_{q>Q}\lambda_q^{2s+2}\|h_q\|_{L^2}^2.
\end{split}
\end{equation}
The estimate of $I_3$ is given by
\begin{equation}\notag
\begin{split}
I_3\leq&\sum_{q>Q}\sum_{p\geq q-2} \lambda_q^{2s} \int_{\mathbb T^n} \left|\Delta_q\left((\nabla\times \widetilde B^{(1)}_p)\times h_p\right)\cdot \nabla\times h_q\right|\, d\vec x\\
\leq&\sum_{q>Q}\sum_{p\geq q-2} \lambda_q^{2s} \|\nabla\times \widetilde B^{(1)}_p\|_{L^r}\|h_p\|_{L^2} \|\nabla\times h_q\|_{L^{\frac{2r}{r-2}}}\\
\lesssim&\sum_{q>Q}\lambda_q^{2s+1+\frac{n}{r}} \|h_q\|_{L^2}\sum_{p\geq q-2}\lambda_p \|\widetilde B^{(1)}_p\|_{L^r}\|h_p\|_{L^2}\\
\leq&\ c_0\mu\sum_{q>Q}\lambda_q^{2s+1+\frac{n}{r}} \|h_q\|_{L^2}\sum_{p\geq q-2}\lambda_p ^{1-\frac{n}{r}}\|h_p\|_{L^2}\\ 
\leq&\ c_0\mu\sum_{q>Q}\lambda_q^{s+1} \|h_q\|_{L^2}\sum_{p\geq q-2}\lambda_p ^{s+1}\|h_p\|_{L^2} \lambda_{q-p} ^{s+\frac{n}{r}}\\
\leq&\ c_0\mu\sum_{q>Q}\lambda_q^{2s+2} \|h_q\|_{L^2}^2
\end{split}
\end{equation}
since $s+\frac{n}{r}>0$. 

We estimate $J$ similarly with the additional help of commutators. Bony's paraproduct decomposition yields
\begin{equation}\notag
\begin{split}
J=&\sum_{q\geq -1}\lambda_q^{2s} \int_{\mathbb T^n} \Delta_q\left((\nabla\times h)\times B^{(2)}\right)\cdot \nabla\times h_q\, d\vec x\\
=&\sum_{q\geq -1}\sum_{|p-q|\leq 2} \lambda_q^{2s} \int_{\mathbb T^n} \Delta_q\left((\nabla\times h_{\leq p-2})\times B^{(2)}_p\right)\cdot \nabla\times h_q\, d\vec x\\
&+\sum_{q\geq -1}\sum_{|p-q|\leq 2} \lambda_q^{2s} \int_{\mathbb T^n} \Delta_q\left((\nabla\times h_{p})\times B^{(2)}_{\leq p-2}\right)\cdot \nabla\times h_q\, d\vec x\\
&+\sum_{q\geq -1}\sum_{p\geq q-2} \lambda_q^{2s} \int_{\mathbb T^n} \Delta_q\left((\nabla\times \widetilde{h}_{p})\times B^{(2)}_p\right)\cdot \nabla\times h_q\, d\vec x\\
=&: J_1+J_2+J_3.
\end{split}
\end{equation}
For $s-\frac{n}{r}<0$, we apply in the order of H\"older's inequality, Bernstein's inequality, (\ref{high-modes}) and Jensen's inequality
\begin{equation}\notag
\begin{split}
|J_1|\leq & \sum_{q>Q}\sum_{|p-q|\leq 2} \lambda_q^{2s} \|\nabla\times h_{\leq p-2}\|_{L^{\frac{2r}{r-2}}}\|B^{(2)}_p\|_{L^r} \|\nabla\times h_q\|_{L^2}\\
\leq &\ c_0\mu\sum_{q>Q}\sum_{|p-q|\leq 2} \lambda_q^{2s+1} \|h_q\|_{L^2}\lambda_p^{-\frac{n}{r}} \sum_{p'\leq p-2} \lambda_{p'}^{1+\frac{n}{r}}\|h_{p'}\|_{L^2}\\
\leq &\ c_0\mu \sum_{q>Q}  \lambda_q^{2s+1-\frac{n}{r}} \|h_q\|_{L^2} \sum_{p'\leq q} \lambda_{p'}^{1+\frac{n}{r}}\|h_{p'}\|_{L^2}\\
\leq &\ c_0\mu \sum_{q>Q}  \lambda_q^{s+1} \|h_q\|_{L^2} \sum_{p'\leq q} \lambda_{p'}^{s+1}\|h_{p'}\|_{L^2} \lambda_{q-p'}^{s-\frac{n}{r}}\\
\leq &\ c_0\mu \sum_{q>Q}  \lambda_q^{2s+2} \|h_q\|_{L^2}^2.
\end{split}
\end{equation}

We use the commutator
\[[\Delta_q, B^{(2)}_{\leq p-2}\times \nabla\times ]h_p=\Delta_q\left(B^{(2)}_{\leq p-2}\times(\nabla\times h_p)\right)- B^{(2)}_{\leq p-2}\times \nabla\times \Delta_q(h_p)\]
to further decompose $J_2$,
%We proved the commutator estimate in \cite{},
%\begin{equation}\notag
%\|[\Delta_q, B^{(2)}_{\leq p-2}\times \nabla\times ]h_p\|_{L^{r_1}}\lesssim \|\nabla\times B^{(2)}_{\leq p-2}\|_{L^{r_2}}\| h_p\|_{L^{r_3}}
%\end{equation}
%with $\frac{1}{r_1}=\frac{1}{r_2}+\frac{1}{r_3}$.
\begin{equation}\notag
\begin{split}
J_2=&-\sum_{q\geq -1}\sum_{|p-q|\leq 2} \lambda_q^{2s} \int_{\mathbb T^n} \Delta_q\left(B^{(2)}_{\leq p-2}\times (\nabla\times h_{p})\right)\cdot \nabla\times h_q\, d\vec x\\
=&-\sum_{q\geq -1}\sum_{|p-q|\leq 2} \lambda_q^{2s} \int_{\mathbb T^n} [\Delta_q, B^{(2)}_{\leq p-2}\times \nabla\times ]h_p \cdot \nabla\times h_q\, d\vec x\\
&-\sum_{q\geq -1}\sum_{|p-q|\leq 2} \lambda_q^{2s} \int_{\mathbb T^n} B^{(2)}_{\leq q-2}\times \nabla\times \Delta_q(h_p) \cdot \nabla\times h_q\, d\vec x\\
&-\sum_{q\geq -1}\sum_{|p-q|\leq 2} \lambda_q^{2s} \int_{\mathbb T^n} \left(B^{(2)}_{\leq p-2}-B^{(2)}_{\leq q-2}\right)\times \nabla\times \Delta_q(h_p) \cdot \nabla\times h_q\, d\vec x\\
=&:-J_{21}-J_{22}-J_{23}.
\end{split}
\end{equation}
Note
\begin{equation}\notag
J_{22}=\sum_{q\geq -1} \lambda_q^{2s} \int_{\mathbb T^n} B^{(2)}_{\leq q-2}\times (\nabla\times h_q) \cdot \nabla\times h_q\, d\vec x
=0.
\end{equation}
For $r> n$ the terms $J_{21}$ and $J_{23}$ are estimated as
\begin{equation}\notag
\begin{split}
|J_{21}|\leq &  \sum_{q>Q} \sum_{\substack {|p-q|\leq 2\\ p>Q+2}}\lambda_q^{2s} \|\nabla\times B^{(2)}_{(Q,p-2]}\|_{L^r} \|h_p\|_{L^{\frac{2r}{r-2}}}\|\nabla\times h_q\|_{L^2}\\
\leq &  \sum_{q>Q} \sum_{\substack {|p-q|\leq 2\\ p>Q+2}}\lambda_q^{2s+1} \|h_q\|_{L^2} \lambda_p^{\frac{n}{r}}  \|h_p\|_{L^2}\sum_{Q<p'\leq p-2} \lambda_{p'} \|B^{(2)}_{(Q,p-2]}\|_{L^r}\\
\leq&\ c_0\mu \sum_{q>Q}\lambda_q^{2s+1+\frac{n}{r}} \|h_q\|_{L^2}^2\sum_{Q<p'\leq q} \lambda_{p'}^{1-\frac{n}{r}}\\
\leq&\ c_0\mu \sum_{q>Q}\lambda_q^{2s+2} \|h_q\|_{L^2}^2\sum_{Q<p'\leq q} \lambda_{p'-q}^{1-\frac{n}{r}}\\
\leq&\ c_0\mu \sum_{q>Q}\lambda_q^{2s+2} \|h_q\|_{L^2}^2,
\end{split}
\end{equation}
\begin{equation}\notag
\begin{split}
|J_{23}|\lesssim &\sum_{q>Q}\lambda_q^{2s+2}\|h_q\|_{L^2} \|h_q\|_{L^{\frac{2r}{r-2}}} \|B^{(2)}_q\|_{L^r}\\
\leq& \ c_0\mu   \sum_{q>Q}\lambda_q^{2s+2}\|h_q\|_{L^2}^2.
\end{split}
\end{equation}
In the end, $J_3$ can be estimated as $I_3$.  Combining the estimates of $I$ and $J$ with (\ref{est-energy1}) we infer for a small enough constant $c_r$
\begin{equation}\notag
\begin{split}
\frac{d}{dt} \sum_{q\geq -1} \lambda_q^{2s}\|h_q\|_{L^2(\mathbb T^n)}^2+\mu \sum_{q\geq -1} \lambda_q^{2s+2}\|h_q\|_{L^2(\mathbb T^n)}^2\leq 0.
\end{split}
\end{equation}
It follows immediately that
\[\lim_{t\to\infty} \|h(t)\|_{H^s}=0, \ \ \ \mbox{for} \ \ -\frac{n}{r}<s<\frac{n}{r}-1.\]
It completes the proof of Theorem \ref{thm-det}.

\bigskip

\section{Low modes regularity condition}
\label{sec-reg}

The proof of Theorem \ref{thm-reg-low} also utilizes the wavenumber splitting approach.  Following standard argument of regularity, it is sufficient to provide a priori estimate for the solution in a subcritical space. In our case, we aim to show the a priori estimate $(a,b)\in H^{s+1}\times H^s$ for some $s>1$. 
%{\color{blue}Question: If $f_1\in L^1(0,T)$, can we show $f_2\in L^1(0,T)$? If so, then low modes condition on the component $b$ is sufficient for regularity. }
Multiplying the equation 
\[a_t=\mu\Delta a-(a_yb_x-a_xb_y)\]
by $\lambda_q^{2s}\Delta a_{qq}$, integrating over $\mathbb T^2$ and taking the summation for $q\geq-1$ gives
\begin{equation}\label{priori-1}
\begin{split}
&\frac{1}{2}\frac{d}{dt} \sum_{q\geq-1} \lambda_q^{2s+2}\int_{\mathbb T^2} a_q^2\, d\vec x+\mu \sum_{q\geq-1} \lambda_q^{2s+4}\int_{\mathbb T^2} a_q^2\, d\vec x\\
=& \sum_{q\geq-1} \lambda_q^{2s}\int_{\mathbb T^2} a_yb_x(a_{xx})_{qq}\, d\vec x+\sum_{q\geq-1} \lambda_q^{2s}\int_{\mathbb T^2} a_yb_x(a_{yy})_{qq}\, d\vec x\\
&- \sum_{q\geq-1} \lambda_q^{2s}\int_{\mathbb T^2} a_xb_y(a_{xx})_{qq}\, d\vec x-\sum_{q\geq-1} \lambda_q^{2s}\int_{\mathbb T^2} a_xb_y(a_{yy})_{qq}\, d\vec x\\
=&: K_1+K_2+K_3+K_4.
\end{split}
\end{equation}
Multiplying the equation 
\[ b_t=\mu\Delta b+(a_y\Delta a_x-a_x\Delta a_y)\]
by $\lambda_q^{2s}b_{qq}$, integrating over $\mathbb T^2$, taking the summation for $q\geq-1$ and applying integration by parts we obtain
\begin{equation}\label{priori-2}
\begin{split}
&\frac{1}{2}\frac{d}{dt} \sum_{q\geq-1} \lambda_q^{2s}\int_{\mathbb T^2} b_q^2\, d\vec x+ \mu\sum_{q\geq-1} \lambda_q^{2s+2}\int_{\mathbb T^2} b_q^2\, d\vec x\\
=& -\sum_{q\geq-1} \lambda_q^{2s}\int_{\mathbb T^2} a_ya_{xx}(b_{x})_{qq}\, d\vec x- \sum_{q\geq-1} \lambda_q^{2s}\int_{\mathbb T^2} a_{xy}a_{xx}b_{qq}\, d\vec x\\
&-\sum_{q\geq-1} \lambda_q^{2s}\int_{\mathbb T^2} a_ya_{yy}(b_{x})_{qq}\, d\vec x-\sum_{q\geq-1} \lambda_q^{2s}\int_{\mathbb T^2} a_{xy}a_{yy}b_{qq}\, d\vec x\\
&+ \sum_{q\geq-1} \lambda_q^{2s}\int_{\mathbb T^2} a_xa_{xx}(b_y)_{qq}\, d\vec x + \sum_{q\geq-1} \lambda_q^{2s}\int_{\mathbb T^2} a_{xy}a_{xx}b_{qq}\, d\vec x\\
&+\sum_{q\geq-1} \lambda_q^{2s}\int_{\mathbb T^2} a_xa_{yy}(b_{y})_{qq}\, d\vec x+\sum_{q\geq-1} \lambda_q^{2s}\int_{\mathbb T^2} a_{xy}a_{yy}b_{qq}\, d\vec x\\
=& -\sum_{q\geq-1} \lambda_q^{2s}\int_{\mathbb T^2} a_ya_{xx}(b_{x})_{qq}\, d\vec x
-\sum_{q\geq-1} \lambda_q^{2s}\int_{\mathbb T^2} a_ya_{yy}(b_{x})_{qq}\, d\vec x\\
&+ \sum_{q\geq-1} \lambda_q^{2s}\int_{\mathbb T^2} a_xa_{xx}(b_y)_{qq}\, d\vec x 
+\sum_{q\geq-1} \lambda_q^{2s}\int_{\mathbb T^2} a_xa_{yy}(b_{y})_{qq}\, d\vec x\\
=&: K_5+K_6+K_7+K_8.
\end{split}
\end{equation}
We will show that there are cancellations in $K_1+K_5$, $K_2+K_6$, $K_3+K_7$ and $K_4+K_8$.    

Applying Bony's paraproduct we obtain
\begin{equation}\notag
\begin{split}
K_1=&\sum_{q\geq -1}\lambda_q^{2s}\int_{\mathbb T^2}(a_yb_x)_q a_{xx,q} \, d\vec x\\
=&\sum_{q\geq -1}\sum_{|p-q|\leq 2}\lambda_q^{2s}\int_{\mathbb T^2}(a_{y,\leq p-2}b_{x,p})_q a_{xx,q} \, d\vec x\\
&+\sum_{q\geq -1}\sum_{|p-q|\leq 2}\lambda_q^{2s}\int_{\mathbb T^2}(b_{x,\leq p-2}a_{y,p})_q a_{xx,q} \, d\vec x\\
&+\sum_{q\geq -1}\sum_{p\geq q- 2}\lambda_q^{2s}\int_{\mathbb T^2}(\widetilde a_{y,p}b_{x,p})_q a_{xx,q} \, d\vec x\\
=&: K_{11}+K_{12}+K_{13}
\end{split}
\end{equation}
and 
\begin{equation}\notag
\begin{split}
K_5=&-\sum_{q\geq-1} \lambda_q^{2s}\int_{\mathbb T^2} (a_ya_{xx})_qb_{x,q}\, d\vec x\\
=&-\sum_{q\geq-1}\sum_{|p-q|\leq 2} \lambda_q^{2s}\int_{\mathbb T^2} (a_{y,\leq p-2}a_{xx, p})_qb_{x,q}\, d\vec x\\
&-\sum_{q\geq-1}\sum_{|p-q|\leq 2} \lambda_q^{2s}\int_{\mathbb T^2} (a_{xx,\leq p-2}a_{y, p})_qb_{x,q}\, d\vec x\\
&-\sum_{q\geq-1}\sum_{p\geq q- 2} \lambda_q^{2s}\int_{\mathbb T^2} (\widetilde a_{y,p}a_{xx, p})_qb_{x,q}\, d\vec x\\
=&: K_{51}+K_{52}+K_{53}.
\end{split}
\end{equation}
Denote the commutators 
\begin{equation}\notag
\begin{split}
[\Delta_q, a_{y,\leq p-2}\partial_x] b_p=&\ \Delta_q\left(a_{y,\leq p-2}b_{x,p}\right)-a_{y,\leq p-2}(b_{x,p})_q,\\
[\Delta_q, a_{y,\leq p-2}\partial_x] a_{x,p}=&\ \Delta_q\left(a_{y,\leq p-2}a_{xx,p}\right)-a_{y,\leq p-2}(a_{xx,p})_q.
\end{split}
\end{equation}
In order to shift derivatives to low modes and explore cancellations, we rewrite $I_{11}$ and $I_{51}$ as
\begin{equation}\notag
\begin{split}
K_{11}=&\sum_{q\geq -1}\sum_{|p-q|\leq 2}\lambda_q^{2s} \int_{\mathbb T^2} [\Delta_q, a_{y,\leq p-2}\partial_x] b_p a_{xx,q} \, d\vec x\\
&+\sum_{q\geq -1}\sum_{|p-q|\leq 2}\lambda_q^{2s} \int_{\mathbb T^2} a_{y,\leq q-2}(b_{x,p})_q a_{xx,q} \, d\vec x\\
&+\sum_{q\geq -1}\sum_{|p-q|\leq 2}\lambda_q^{2s} \int_{\mathbb T^2} \left(a_{y,\leq q-2}-a_{y,\leq p-2}\right)(b_{x,p})_q a_{xx,q} \, d\vec x\\
=&:K_{111}+K_{112}+K_{113},
\end{split}
\end{equation}
\begin{equation}\notag
\begin{split}
K_{51}=&-\sum_{q\geq -1}\sum_{|p-q|\leq 2}\lambda_q^{2s} \int_{\mathbb T^2} [\Delta_q, a_{y,\leq p-2}\partial_x] a_{x,p} b_{x,q} \, d\vec x\\
&-\sum_{q\geq -1}\sum_{|p-q|\leq 2}\lambda_q^{2s} \int_{\mathbb T^2} a_{y,\leq q-2}(a_{xx,p})_q b_{x,q} \, d\vec x\\
&-\sum_{q\geq -1}\sum_{|p-q|\leq 2}\lambda_q^{2s} \int_{\mathbb T^2} \left(a_{y,\leq q-2}-a_{y,\leq p-2}\right)(a_{xx,p})_q b_{x,q} \, d\vec x\\
=&:K_{511}+K_{512}+K_{513}.
\end{split}
\end{equation}
In observation of $\sum_{|p-q|\leq 2}(b_{x,p})_q=b_{x,q}$ and $\sum_{|p-q|\leq 2}(a_{xx,p})_q=a_{xx,q}$, we have the cancellation 
\[K_{112}+K_{512}=0.\]
The other terms in $K_{11}+K_{51}$ are estimated as follows. Splitting the frequency into high and low parts in $K_{111}$ yields
\begin{equation}\notag
\begin{split}
K_{111}=&\sum_{p>Q}\sum_{|q-p|\leq 2}\lambda_q^{2s} \int_{\mathbb T^2} [\Delta_q, a_{y,\leq p-2}\partial_x] b_p a_{xx,q} \, d\vec x\\
&+\sum_{-1\leq p\leq Q}\sum_{|q-p|\leq 2}\lambda_q^{2s} \int_{\mathbb T^2} [\Delta_q, a_{y,\leq p-2}\partial_x] b_p a_{xx,q} \, d\vec x\\
=&:K_{1111}+K_{1112}.
\end{split}
\end{equation}
Applying the commutator estimate in Lemma \ref{le-comm2}, the definition of the wavenumber $\lambda_{Q(b)}$, and Bernstein's inequality we obtain
\begin{equation}\notag
\begin{split}
|K_{1111}|\leq & \sum_{p>Q}\sum_{|q-p|\leq 2}\lambda_q^{2s}\|a_{xxy,\leq p-2}\|_{L^{\frac{2r}{r-2}}}\| b_p\|_{L^r} \|a_{x,q}\|_{L^2}\\
\lesssim &  \sum_{p>Q}\lambda_p^{2s+1}\|a_{p}\|_{L^2}\| b_p\|_{L^r}\sum_{p'\leq p-2} \lambda_{p'}^{3+\frac{2}{r}}\|a_{p'}\|_{L^2}\\
\leq &\ c_r\mu \sum_{p>Q}\lambda_p^{2s+1-\frac{2}{r}}\|a_{p}\|_{L^2}\sum_{p'\leq p-2} \lambda_{p'}^{3+\frac{2}{r}}\|a_{p'}\|_{L^2}.
\end{split}
\end{equation}
For $1<s<1+\frac{2}{r}$, we continue to use Jensen's inequality
\begin{equation}\notag
\begin{split}
|K_{1111}|\leq &\ c_r\mu \sum_{p>Q}\lambda_p^{s+2}\|a_{p}\|_{L^2}\sum_{p'\leq p-2} \lambda_{p'}^{s+2}\|a_{p'}\|_{L^2}\lambda_{p-p'}^{s-1-\frac{2}{r}}\\
\leq &\ c_r\mu \sum_{p>Q}\lambda_p^{2s+4}\|a_{p}\|_{L^2}^2+c_r\mu \sum_{p>Q} \sum_{p'\leq p-2} \lambda_{p'}^{2s+4}\|a_{p'}\|_{L^2}^2\lambda_{p-p'}^{s-1-\frac{2}{r}}\\
\leq &\ c_r\mu \sum_{p>Q}\lambda_p^{2s+4}\|a_{p}\|_{L^2}^2+c_r\mu \sum_{p'\geq -1} \lambda_{p'}^{2s+4}\|a_{p'}\|_{L^2}^2 \sum_{p\geq p'+2} \lambda_{p-p'}^{s-1-\frac{2}{r}}\\
\leq &\ c_r\mu \sum_{q\geq -1}\lambda_p^{2s+4}\|a_{q}\|_{L^2}^2. 
\end{split}
\end{equation}
On the other hand, it follows from the commutator estimate, the definition of $f(t)$ and Jensen's inequality that
\begin{equation}\notag
\begin{split}
|K_{1112}|\leq & \sum_{-1\leq p\leq Q}\sum_{|q-p|\leq 2}\lambda_q^{2s}\|a_{xxy,\leq p-2}\|_{L^2}\| b_p\|_{L^\infty} \|a_{x,q}\|_{L^2}\\
\lesssim& \sum_{-1\leq p\leq Q}\lambda_p^{2s+1}\|a_{p}\|_{L^2} \| b_p\|_{L^\infty}\sum_{p'\leq p-2}\lambda_{p'}^3\|a_{p'}\|_{L^2}\\
\leq &\ Cf_1(t) \sum_{-1\leq p\leq Q}\lambda_p^{2s-1}\|a_{p}\|_{L^2} \sum_{p'\leq p-2}\lambda_{p'}^3\|a_{p'}\|_{L^2}\\
\leq &\ Cf_1(t) \sum_{-1\leq p\leq Q}\lambda_p^{s+1}\|a_{p}\|_{L^2} \sum_{p'\leq p-2}\lambda_{p'}^{s+1}\|a_{p'}\|_{L^2}\lambda_{p-p'}^{s-2}\\
\leq &\ Cf_1(t) \sum_{-1\leq p\leq Q}\lambda_p^{2s+2}\|a_{p}\|_{L^2}^2
\end{split}
\end{equation}
for $1<s<2$. 

The frequency interaction of the term $K_{113}$ is weaker than $K_{111}$, and hence the estimate of it is simpler. Splitting the frequency in $K_{113}$ gives rise to 
\begin{equation}\notag
\begin{split}
K_{113}=& \sum_{p\geq -1}\sum_{|p-q|\leq 2}\lambda_q^{2s} \int_{\mathbb T^2} \left(a_{y,\leq q-2}-a_{y,\leq p-2}\right)(b_{x,p})_q a_{xx,q} \, d\vec x\\
=&\sum_{p>Q}\sum_{|p-q|\leq 2}\lambda_q^{2s} \int_{\mathbb T^2} \left(a_{y,\leq q-2}-a_{y,\leq p-2}\right)(b_{x,p})_q a_{xx,q} \, d\vec x\\
&+\sum_{-1\leq p\leq Q}\sum_{|p-q|\leq 2}\lambda_q^{2s} \int_{\mathbb T^2} \left(a_{y,\leq q-2}-a_{y,\leq p-2}\right)(b_{x,p})_q a_{xx,q} \, d\vec x\\
=&: I_{1131}+I_{1132},
\end{split}
\end{equation}
while
\begin{equation}\notag
\begin{split}
|K_{1131}|\leq & \sum_{p>Q}\sum_{|q-p|\leq 2}\lambda_q^{2s} \|a_{y,q}\|_{L^2}\|b_{x,p}\|_{L^r}\|a_{xx,q}\|_{L^{\frac{2r}{r-2}}}\\
\lesssim& \sum_{p>Q}\sum_{|q-p|\leq 2}\lambda_q^{2s+3+\frac{2}{r}} \|a_q\|_{L^2}^2\|b_{x,p}\|_{L^r}\\
\leq&\ c_r\mu \sum_{q>Q-2}\lambda_q^{2s+4} \|a_q\|_{L^2}^2\\
\end{split}
\end{equation}
and 
\begin{equation}\notag
\begin{split}
|K_{1132}|\leq & \sum_{p\leq Q}\sum_{|q-p|\leq 2}\lambda_q^{2s} \|a_{y,q}\|_{L^2}\|b_{x,p}\|_{L^\infty}\|a_{xx,q}\|_{L^2}\\
\leq&\ Cf_1(t) \sum_{q>Q-2}\lambda_q^{2s+2} \|a_q\|_{L^2}^2.
\end{split}
\end{equation}
The estimate of $K_{511}$ follows from the frequency splitting
\begin{equation}\notag
\begin{split}
K_{511}=&-\sum_{q>Q}\sum_{|p-q|\leq 2}\lambda_q^{2s} \int_{\mathbb T^2} [\Delta_q, a_{y,\leq p-2}\partial_x] a_{x,p} b_{x,q} \, d\vec x\\
&-\sum_{q\leq Q}\sum_{|p-q|\leq 2}\lambda_q^{2s} \int_{\mathbb T^2} [\Delta_q, a_{y,\leq p-2}\partial_x] a_{x,p} b_{x,q} \, d\vec x\\
=&:K_{5111}+K_{5112},
\end{split}
\end{equation}
and commutator estimate in Lemma \ref{le-comm2}, Bernstein's inequality and Jensen's inequality
\begin{equation}\notag
\begin{split}
|K_{5111}|\leq & \sum_{q>Q}\sum_{|p-q|\leq 2}\lambda_q^{2s}\|a_{xxy,\leq p-2}\|_{L^{\frac{2r}{r-2}}} \|a_{x,p}\|_{L^2}\|b_{q}\|_{L^r}\\
\lesssim & \sum_{q>Q}\sum_{|p-q|\leq 2}\lambda_q^{2s+1} \|a_{p}\|_{L^2}\|b_{q}\|_{L^r}\sum_{p'\leq p-2}\lambda_{p'}^{3+\frac{2}{r}} \|a_{p'}\|_{L^2} \\
\leq & \ c_r\mu\sum_{q>Q}\sum_{|p-q|\leq 2}\lambda_q^{2s+1-\frac2r} \|a_{p}\|_{L^2}\sum_{p'\leq p-2}\lambda_{p'}^{3+\frac{2}{r}} \|a_{p'}\|_{L^2} \\
\leq & \ c_r\mu\sum_{q>Q}\lambda_q^{s+2} \|a_{q}\|_{L^2}\sum_{p'\leq q}\lambda_{p'}^{s+2} \|a_{p'}\|_{L^2}\lambda_{q-p'}^{s-1-\frac{2}{r}}  \\
\leq & \ c_r\mu\sum_{q\geq -1}\lambda_q^{2s+4} \|a_{q}\|_{L^2}^2
\end{split}
\end{equation}
for $1<s<1+\frac2r$, and 
\begin{equation}\notag
\begin{split}
|K_{5112}|\leq & \sum_{q\leq Q}\sum_{|p-q|\leq 2}\lambda_q^{2s}\|a_{xxy,\leq p-2}\|_{L^\infty} \|a_{x,p}\|_{L^2}\|b_{q}\|_{L^2}\\
\lesssim & \sum_{q\leq Q}\sum_{|p-q|\leq 2}\lambda_q^{2s+1} \|a_{p}\|_{L^2}\|b_{q}\|_{L^2}\sum_{p'\leq p-2}\lambda_{p'}^{3} \|a_{p'}\|_{L^2} \\
\leq & \ C f_1(t)\sum_{q\leq Q}\sum_{|p-q|\leq 2}\lambda_q^{2s-1} \|a_{p}\|_{L^2}\sum_{p'\leq p-2}\lambda_{p'}^{3} \|a_{p'}\|_{L^2} \\
\leq & \ Cf_1(t)\sum_{q\leq Q+2}\lambda_q^{s+1} \|a_{q}\|_{L^2}\sum_{p'\leq q}\lambda_{p'}^{s+1} \|a_{p'}\|_{L^2}\lambda_{q-p'}^{s-2}  \\
\leq & \ Cf_1(t)\sum_{q\leq Q+2}\lambda_q^{2s+2} \|a_{q}\|_{L^2}^2
\end{split}
\end{equation}
for $1<s<2$. The estimate of $K_{513}$ is similar to that of $K_{113}$. 

Therefore we conclude from the analysis above that
\begin{equation}\notag
\begin{split}
|K_{11}+K_{51}|\leq c_r\mu\sum_{q\geq -1}\lambda_q^{2s+4} \|a_{q}\|_{L^2}^2+Cf_1(t)\sum_{q\leq Q}\lambda_q^{2s+2} \|a_{q}\|_{L^2}^2.
\end{split}
\end{equation}
The estimate of $K_{12}+K_{52}$ is similar to that of $K_{11}+K_{51}$. The difference is that the frequency splitting acts on the term $a_{y,p}$; hence the wavenumber in (\ref{wave-a}) is used to separate the high and low frequency parts. We obtain the estimate 
\begin{equation}\notag
\begin{split}
|K_{12}+K_{52}|\leq c_r\mu\sum_{q\geq -1}\lambda_q^{2s+4} \|a_{q}\|_{L^2}^2+Cf_2(t)\sum_{q\leq Q}\lambda_q^{2s} \|b_{q}\|_{L^2}^2.
\end{split}
\end{equation}

The estimate of $K_{13}$ is as follows, with frequency splitting first 
\begin{equation}\notag
\begin{split}
K_{13}=&\sum_{p\geq -1}\sum_{q\leq p+2}\lambda_q^{2s}\int_{\mathbb T^2}(\widetilde a_{y,p}b_{x,p})_q a_{xx,q} \, d\vec x\\
=& \sum_{p>Q}\sum_{q\leq p+2}\lambda_q^{2s}\int_{\mathbb T^2}(\widetilde a_{y,p}b_{x,p})_q a_{xx,q} \, d\vec x\\
&+\sum_{p\leq Q}\sum_{q\leq p+2}\lambda_q^{2s}\int_{\mathbb T^2}(\widetilde a_{y,p}b_{x,p})_q a_{xx,q} \, d\vec x\\
=&: K_{131}+K_{132}
\end{split}
\end{equation}
and then
\begin{equation}\notag
\begin{split}
|K_{131}|\leq & \sum_{p>Q}\sum_{q\leq p+2}\lambda_q^{2s}\|\widetilde a_{y,p}\|_{L^2}\|b_{x,p}\|_{L^r}\|a_{xx,q}\|_{L^{\frac{2r}{r-2}}}\\
\lesssim& \sum_{p>Q}\lambda_p^2\|\widetilde a_{p}\|_{L^2}\|b_{p}\|_{L^r}\sum_{q\leq p+2}\lambda_q^{2s+2+\frac2r}\|a_{q}\|_{L^2}\\
\leq&\ c_r\mu \sum_{p>Q}\lambda_p^{2-\frac2r}\|\widetilde a_{p}\|_{L^2}\sum_{q\leq p+2}\lambda_q^{2s+2+\frac2r}\|a_{q}\|_{L^2}\\
\leq&\ c_r\mu \sum_{p>Q}\lambda_p^{s+2}\|\widetilde a_{p}\|_{L^2}\sum_{q\leq p+2}\lambda_q^{s+2}\|a_{q}\|_{L^2}\lambda_{q-p}^{s+\frac2r}\\
\leq&\ c_r\mu \sum_{p\geq -1}\lambda_p^{2s+4}\|a_{p}\|_{L^2}^2
\end{split}
\end{equation}
since $s+\frac2r>0$, and 
\begin{equation}\notag
\begin{split}
|K_{132}|\leq & \sum_{p\leq Q}\sum_{q\leq p+2}\lambda_q^{2s}\|\widetilde a_{y,p}\|_{L^2}\|b_{x,p}\|_{L^\infty}\|a_{xx,q}\|_{L^2}\\
\lesssim& \sum_{p\leq Q}\lambda_p^2\|\widetilde a_{p}\|_{L^2}\|b_{p}\|_{L^\infty}\sum_{q\leq p+2}\lambda_q^{2s+2}\|a_{q}\|_{L^2}\\
\leq&\ Cf_1(t) \sum_{p\leq Q}\|\widetilde a_{p}\|_{L^2}\sum_{q\leq p+2}\lambda_q^{2s+2}\|a_{q}\|_{L^2}\\
\leq&\ Cf_1(t) \sum_{p\leq Q}\lambda_p^{s+1}\|\widetilde a_{p}\|_{L^2}\sum_{q\leq p+2}\lambda_q^{s+1}\|a_{q}\|_{L^2}\lambda_{q-p}^{s+1}\\
\leq&\ Cf_1(t) \sum_{p\leq Q+2}\lambda_p^{2s+2}\|a_{p}\|_{L^2}^2
\end{split}
\end{equation}
for $s+1>0$. The estimate of $K_{53}$ is similar to that of $K_{13}$. Therefore combining the analysis above we have for $1<s<1+\frac2r$
\begin{equation}\label{priori-3}
\begin{split}
|K_{1}+K_{5}|\leq &\ Cc_r\mu\sum_{q\geq -1}\lambda_q^{2s+4} \|a_{q}\|_{L^2}^2\\
&+C\left(f_1(t)+f_2(t)\right)\sum_{q\geq -1}\left(\lambda_q^{2s+2} \|a_{q}\|_{L^2}^2+\lambda_q^{2s} \|b_{q}\|_{L^2}^2\right).
\end{split}
\end{equation}

In the end, we observe that $K_2+K_6$, $K_3+K_7$ and $K_4+K_8$ can be estimated analogously as for $K_1+K_5$, by exploring the cancellations,
\begin{equation}\label{priori-4}
\begin{split}
&|K_{2}+K_{6}|+|K_{3}+K_{7}|+|K_{4}+K_{8}|\\
\leq &\ Cc_r\mu\sum_{q\geq -1}\lambda_q^{2s+4} \|a_{q}\|_{L^2}^2\\
&+C\left(f_1(t)+f_2(t)\right)\sum_{q\geq -1}\left(\lambda_q^{2s+2} \|a_{q}\|_{L^2}^2+\lambda_q^{2s} \|b_{q}\|_{L^2}^2\right).
\end{split}
\end{equation}
It follows from (\ref{priori-1}),(\ref{priori-2}) and (\ref{priori-3}) that, for any $r\in[2,\infty)$ and a small enough constant $c_r$
\begin{equation}\notag
\begin{split}
&\frac{d}{dt} \sum_{q\geq -1}\left(\lambda_q^{2s+2} \|a_{q}\|_{L^2}^2+\lambda_q^{2s} \|b_{q}\|_{L^2}^2\right)+\mu\sum_{q\geq -1}\left(\lambda_q^{2s+4} \|a_{q}\|_{L^2}^2+\lambda_q^{2s+2} \|b_{q}\|_{L^2}^2\right)\\
\leq &\ C\left(f_1(t)+f_2(t)\right)\sum_{q\geq -1}\left(\lambda_q^{2s+2} \|a_{q}\|_{L^2}^2+\lambda_q^{2s} \|b_{q}\|_{L^2}^2\right).
\end{split}
\end{equation}
Thus in view of the assumption $f_1,f_2\in L^1(0,T)$ and Gr\"onwall's inequality we obtain 
\[ (a(t),b(t))\in H^{s+1}\times H^s, \ \ \forall \ \ t\in[0,T]\]
for some $s>1$. This concludes the proof.

\bigskip

\section{Proof of Theorem \ref{thm-criterion}}
\label{sec-criterion}

%It follows from standard analysis that the Ladyzhenskaya-Prodi-Serrin type of regularity criterion holds for the electron MHD (\ref{emhd}). In particular, for a solution $B=B(x,y,t)$ of system (\ref{emhd}), if 
%\begin{equation}\label{LPS}
%\int_0^T\|\nabla B\|_{L^p(\mathbb T^2)}^q\, dt<\infty \ \ \ \mbox{with} \ \ \ \frac{2}{q}+\frac{2}{p}\leq 1, \ \ 2<p\leq \infty,
%\end{equation}
%the solution is regular on $[0,T]$. Note that if $a=a(x,y,t)$ and $b=b(x,y,t)$ satisfy 
%\begin{equation}\label{priori}
%\nabla(\nabla\times (a\vec e_z))\in L^2(0,T; L^\infty(\mathbb T^2)), \ \ \ \nabla(b\vec e_z)\in L^2(0,T; L^\infty(\mathbb T^2)),
%\end{equation}
%the vector field $B=\nabla\times (a\vec e_z)+b\vec e_z$ satisfies (\ref{LPS}) with $p=\infty$ and $q=2$. Therefore in order to prove Theorem \ref{thm-criterion} under the assumption (\ref{criterion}), it is sufficient to show the a priori estimate for the function $a$ in (\ref{priori}).
In order to prove Theorem \ref{thm-criterion}, we adapt the standard energy method rather than the harmonic analysis techniques used in previous sections. An iteration process will be taken to establish higher order energy estimates. 
 Acting gradient $\nabla$ on the equation of $a$ in (\ref{emhd1}) and taking dot product with $|\nabla a|^{p-2}\nabla a$ for any $p\geq 2$, it then follows from integrating over $\mathbb T^2$ that
\begin{equation}\label{est-ap1}
\begin{split}
&\frac{1}{p}\frac{d}{dt}\int_{\mathbb T^2} |\nabla a|^p\, d\vec x+\mu (p-1) \int_{\mathbb T^2} |\nabla \nabla a|^2 |\nabla a|^{p-2}\, d\vec x\\
=&- \int_{\mathbb T^2}\nabla(a_yb_x)\cdot \nabla a |\nabla a|^{p-2}\, d\vec x+\int_{\mathbb T^2}\nabla(a_xb_y)\cdot \nabla a |\nabla a|^{p-2}\, d\vec x.
\end{split}
\end{equation}
Invoking integration by parts in the first integral on the right hand side of (\ref{est-ap1}) gives
\begin{equation}\notag
- \int_{\mathbb T^2}\nabla(a_yb_x)\cdot \nabla a |\nabla a|^{p-2}\, d\vec x= \int_{\mathbb T^2}a_yb_x\nabla\cdot\left(\nabla a |\nabla a|^{p-2}\right)\, d\vec x. 
\end{equation}
Hence we have the estimate by using Young's inequality
\begin{equation}\notag
\begin{split}
&\left| \int_{\mathbb T^2}\nabla(a_yb_x)\cdot \nabla a |\nabla a|^{p-2}\, d\vec x \right|\\
=&\left| \int_{\mathbb T^2}a_yb_x\nabla\cdot\left(\nabla a |\nabla a|^{p-2}\right)\, d\vec x \right|\\
\leq&\ (p-1) \int_{\mathbb T^2}|b_x| |\nabla a|^{p-1} |\nabla\nabla a|\, d\vec x\\
\leq&\ \frac{1}{4} (p-1)\int_{\mathbb T^2} |\nabla \nabla a|^2 |\nabla a|^{p-2}\, d\vec x+(p-1) \int_{\mathbb T^2}|b_x|^2 |\nabla a|^{p}\, d\vec x.
\end{split}
\end{equation}
Analogously the second integral on the right hand side of (\ref{est-ap1}) has the estimate 
\begin{equation}\notag
\begin{split}
&\left| \int_{\mathbb T^2}\nabla(a_xb_y)\cdot \nabla a |\nabla a|^{p-2}\, d\vec x \right|\\
\leq&\ \frac{1}{4} (p-1)\int_{\mathbb T^2} |\nabla \nabla a|^2 |\nabla a|^{p-2}\, d\vec x+(p-1) \int_{\mathbb T^2}|b_y|^2 |\nabla a|^{p}\, d\vec x.
\end{split}
\end{equation}
Combining the estimates above with (\ref{est-ap1}) yields
\begin{equation}\label{est-ap2}
\begin{split}
&\frac{d}{dt}\int_{\mathbb T^2} |\nabla a|^p\, d\vec x+\frac12\mu p(p-1) \int_{\mathbb T^2} |\nabla \nabla a|^2 |\nabla a|^{p-2}\, d\vec x\\
\leq &\ 2p(p-1) \|\nabla b\|_{L^\infty}^2\int_{\mathbb T^2} |\nabla a|^{p}\, d\vec x.
\end{split}
\end{equation}
It follows from Gr\"onwall's inequality, (\ref{est-ap2}) and the assumption (\ref{criterion}) that
\begin{equation}\label{est-ap3}
\nabla a\in L^\infty(0,T; L^p(\mathbb T^2)), \ \ 2\leq p\leq \infty
\end{equation}
for any time $T>0$. 

We proceed to estimate the higher order derivatives $\Delta a$ and $\nabla b$. A standard energy argument for the system (\ref{emhd1}) leads to 
\begin{equation}\label{est-ab1}
\begin{split}
&\frac12\frac{d}{dt} \int_{\mathbb T^2}|\Delta a|^2+|\nabla b|^2\, d\vec x+\mu \int_{\mathbb T^2}|\nabla \Delta a|^2+|\Delta b|^2\, d\vec x\\
=&- \int_{\mathbb T^2}\Delta(a_yb_x)\Delta a\, d\vec x+\int_{\mathbb T^2}\Delta(a_xb_y)\Delta a\, d\vec x\\
&- \int_{\mathbb T^2}a_y\Delta a_x\Delta b\, d\vec x+\int_{\mathbb T^2}a_x\Delta a_y\Delta b\, d\vec x\\
=:&\ L_1+L_2+L_3+L_4.
\end{split}
\end{equation}
%Note 
%\[\Delta(a_yb_x)=\Delta a_y b_x+a_y\Delta b_x+2a_{xy} b_{xx}+2a_{yy}b_{xy},\]
%\[\Delta(a_xb_y)=\Delta a_x b_y+a_x\Delta b_y+2a_{xx} b_{xy}+2a_{xy}b_{yy}.\]
Applying integration by parts $L_1$ can be written as
\begin{equation}\notag
\begin{split}
L_1=&-\int_{\mathbb T^2}(a_yb_x)_{xx}\Delta a\, d\vec x-\int_{\mathbb T^2}(a_yb_x)_{yy}\Delta a\, d\vec x\\
=& \int_{\mathbb T^2}(a_yb_x)_{x}\Delta a_x\, d\vec x+\int_{\mathbb T^2}(a_yb_x)_{y}\Delta a_y\, d\vec x\\
=&  \int_{\mathbb T^2}a_yb_{xx}\Delta a_x\, d\vec x+\int_{\mathbb T^2}a_{xy}b_{x}\Delta a_x\, d\vec x\\
&+\int_{\mathbb T^2}a_{yy}b_{x}\Delta a_y\, d\vec x+\int_{\mathbb T^2}a_{y}b_{xy}\Delta a_y\, d\vec x\\
=:& \ L_{11}+L_{12}+L_{13}+L_{14},
\end{split}
\end{equation}
and similarly
\begin{equation}\notag
\begin{split}
L_2=&  -\int_{\mathbb T^2}a_xb_{yy}\Delta a_y\, d\vec x-\int_{\mathbb T^2}a_{xx}b_{y}\Delta a_x\, d\vec x\\
&-\int_{\mathbb T^2}a_{xy}b_{x}\Delta a_y\, d\vec x-\int_{\mathbb T^2}a_{x}b_{xy}\Delta a_x\, d\vec x\\
=:& \ L_{21}+L_{22}+L_{23}+L_{24}.
\end{split}
\end{equation}
On the other hand we see
\begin{equation}\notag
\begin{split}
L_3=&-\int_{\mathbb T^2}a_y\Delta a_x b_{xx}\, d\vec x-\int_{\mathbb T^2}a_y\Delta a_x b_{yy}\, d\vec x=: L_{31}+L_{32},\\
L_4=&\int_{\mathbb T^2}a_x\Delta a_y b_{xx}\, d\vec x+\int_{\mathbb T^2}a_x\Delta a_y b_{yy}\, d\vec x=: L_{41}+L_{42}.
\end{split}
\end{equation}
We observe some cancellations among $L_1+L_3$ and $L_2+L_4$, i.e.
\[L_{11}+L_{31}=0, \ \ \ L_{21}+L_{42}=0.\]
We then apply integration by parts to $L_{14}$, $L_{24}$, $L_{32}$ and $L_{41}$ in order to get rid of second order derivatives of the function $b$, 
\begin{equation}\notag
\begin{split}
L_{14}+L_{32}=&-\int_{\mathbb T^2}a_{xy}b_{y}\Delta a_y\, d\vec x-\int_{\mathbb T^2}a_{y}b_{y}\Delta a_{xy}\, d\vec x\\
&+\int_{\mathbb T^2}a_{yy}b_{y}\Delta a_x\, d\vec x+\int_{\mathbb T^2}a_{y}b_{y}\Delta a_{xy}\, d\vec x\\
=&-\int_{\mathbb T^2}a_{xy}b_{y}\Delta a_y\, d\vec x+\int_{\mathbb T^2}a_{yy}b_{y}\Delta a_x\, d\vec x,
\end{split}
\end{equation}
\begin{equation}\notag
\begin{split}
L_{24}+L_{41}=&\int_{\mathbb T^2}a_{xy}b_{x}\Delta a_x\, d\vec x+\int_{\mathbb T^2}a_{x}b_{x}\Delta a_{xy}\, d\vec x\\
&-\int_{\mathbb T^2}a_{xx}b_{x}\Delta a_y\, d\vec x-\int_{\mathbb T^2}a_{x}b_{x}\Delta a_{xy}\, d\vec x\\
=&\int_{\mathbb T^2}a_{xy}b_{x}\Delta a_x\, d\vec x-\int_{\mathbb T^2}a_{xx}b_{x}\Delta a_y\, d\vec x,
\end{split}
\end{equation}
where more cancellations were explored. Putting the analysis above together we obtain 
\begin{equation}\label{est-ab2}
\begin{split}
L_1+L_2+L_3+L_4=& \ L_{12}+L_{13}+L_{22}+L_{23}\\
&-\int_{\mathbb T^2}a_{xy}b_{y}\Delta a_y\, d\vec x+\int_{\mathbb T^2}a_{yy}b_{y}\Delta a_x\, d\vec x\\
&+\int_{\mathbb T^2}a_{xy}b_{x}\Delta a_x\, d\vec x-\int_{\mathbb T^2}a_{xx}b_{x}\Delta a_y\, d\vec x.
\end{split}
\end{equation}
Note that all of the integrals on the right hand side of (\ref{est-ab2}) have the same order distributions of derivatives among the triple terms. Thus we only show the estimate for $L_{12}$ and the other terms can be handled analogously. It follows from H\"older's inequality that for any $r>2$
\begin{equation}\notag
|L_{12}|=\left| \int_{\mathbb T^2}a_{xy}b_{x}\Delta a_x\, d\vec x \right|
\leq \|\Delta a_x\|_{L^2}\|b_x\|_{L^r}\|a_{xy}\|_{L^{\frac{2r}{r-2}}}.
\end{equation}
Applying Gagliardo-Nirenberg's interpolation inequality 
and the fact
\[\|a_{xy}\|_{L^2}^2=\int_{\mathbb T^2} a_{xy}a_{xy}\, d\vec x= \int_{\mathbb T^2} a_{xx}a_{yy}\, d\vec x\leq \frac12\left(\|a_{xx}\|_{L^2}^2+\|a_{yy}\|_{L^2}^2\right)\]
we deduce
\begin{equation}\notag
\|a_{xy}\|_{L^{\frac{2r}{r-2}}}\leq C \|a_{xy}\|_{L^2}^{\frac{r-2}{r}}\|\nabla a_{xy}\|_{L^2}^{\frac{2}{r}}
\leq C \|\Delta a\|_{L^2}^{\frac{r-2}{r}}\|\nabla \Delta a\|_{L^2}^{\frac{2}{r}}
\end{equation}
for some constant $C>0$. Therefore we conclude by using Young's inequality
\begin{equation}\notag
\begin{split}
|L_{12}|\leq &\ C \|\nabla \Delta a\|_{L^2}^{\frac{r+2}{r}}\| \Delta a\|_{L^2}^{\frac{r-2}{r}}\|\nabla b\|_{L^r}\\
\leq&\ \frac{1}{16}\|\nabla\Delta a\|_{L^2}^2+C\|\nabla b\|_{L^r}^{\frac{2r}{r-2}}\|\Delta a\|_{L^2}^2.
\end{split}
\end{equation}
Hence we have
\begin{equation}\notag
\left| L_1+L_2+L_3+L_4 \right|\leq \frac{1}{2}\|\nabla\Delta a\|_{L^2}^2+C\|\nabla b\|_{L^r}^{\frac{2r}{r-2}}\|\Delta a\|_{L^2}^2
\end{equation}
which together with (\ref{est-ab1}) implies
\begin{equation}\label{est-ab3}
\begin{split}
\frac{d}{dt} \int_{\mathbb T^2}|\Delta a|^2+|\nabla b|^2\, d\vec x+ \mu\int_{\mathbb T^2}|\nabla \Delta a|^2+|\Delta b|^2\, d\vec x\leq C\|\nabla b\|_{L^r}^{\frac{2r}{r-2}}\|\Delta a\|_{L^2}^2.
\end{split}
\end{equation}
Invoking Gr\"onwall's inequality to (\ref{est-ab3}) we claim for any $T>0$
\begin{equation}\label{est-ab5}
\Delta a, \nabla b\in L^\infty(0,T; L^2(\mathbb T^2))\cap L^2(0,T; H^1(\mathbb T^2))
\end{equation}
provided that 
\[\|\nabla b\|_{L^r}\in L^{\frac{2r}{r-2}}(0,T) \]
or more generally 
\[\nabla b\in L^s(0,T; L^r(\mathbb T^2)) \ \ \ \mbox{with} \ \ \frac{2}{s}+\frac{2}{r}\leq 1 \ \ \forall \ \ r\in(2,\infty].\]

%It looks like we need to estimate 
%\begin{equation}\label{est-ab4}
%\frac{d}{dt} \int_{\mathbb T^2}|a_{xy}|^2 \, d\vec x
%\end{equation}
%If that can be estimated as expected, we obtain
%\begin{equation}\label{est-ab5}
%\begin{split}
%&\nabla \nabla a\in L^\infty(0,T; L^2(\mathbb T^2))\cap L^2(0,T; H^1(\mathbb T^2)),\\
%&\nabla b\in L^\infty(0,T; L^2(\mathbb T^2))\cap L^2(0,T; H^1(\mathbb T^2)).
%\end{split}
%\end{equation}
%{\color{blue} Is the estimate $\nabla \nabla a\in L^2(0,T; H^1(\mathbb T^2))$ enough to show $\nabla \nabla a\in L^2(0,T; L^\infty(\mathbb T^2))$?}

We continue to the estimates of the higher order derivatives $\nabla \Delta a$ and $\Delta b$, starting with
\begin{equation}\label{est-ab8}
\begin{split}
&\frac12\frac{d}{dt} \int_{\mathbb T^2}|\nabla \Delta a|^2+|\Delta b|^2\, d\vec x+\mu \int_{\mathbb T^2}|\Delta \Delta a|^2+|\nabla\Delta b|^2\, d\vec x\\
=& \int_{\mathbb T^2}\Delta(a_yb_x)\Delta\Delta a\, d\vec x-\int_{\mathbb T^2}\Delta(a_xb_y)\Delta\Delta a\, d\vec x\\
&+ \int_{\mathbb T^2}a_y\Delta a_x\Delta\Delta b\, d\vec x-\int_{\mathbb T^2}a_x\Delta a_y\Delta\Delta b\, d\vec x\\
=:&\ M_1+M_2+M_3+M_4.
\end{split}
\end{equation}
 Note 
\begin{equation}\notag
\begin{split}
\Delta(a_yb_x)=&\ \Delta a_y b_x+a_y\Delta b_x+2a_{xy} b_{xx}+2a_{yy}b_{xy},\\
\Delta(a_xb_y)=&\ \Delta a_x b_y+a_x\Delta b_y+2a_{xx} b_{xy}+2a_{xy}b_{yy},\\
\Delta(a_y\Delta a_x)=&\ \Delta a_y \Delta a_x+a_y\Delta \Delta a_x+2a_{xy} \Delta a_{xx}+2a_{yy}\Delta a_{xy},\\
\Delta(a_x\Delta a_y)=&\ \Delta a_x \Delta a_y+a_x\Delta \Delta a_y+2a_{xx} \Delta a_{xy}+2a_{xy}\Delta a_{yy}.\\
\end{split}
\end{equation}
Applying integration by parts we can rewrite
\begin{equation}\notag
\begin{split}
M_1=&\ \int_{\mathbb T^2}\Delta a_yb_x\Delta\Delta a\, d\vec x+ \int_{\mathbb T^2}a_y\Delta b_x\Delta\Delta a\, d\vec x\\
&+2 \int_{\mathbb T^2}a_{xy}b_{xx}\Delta\Delta a\, d\vec x+2 \int_{\mathbb T^2}a_{yy}b_{xy}\Delta\Delta a\, d\vec x\\
=&\ \int_{\mathbb T^2}\Delta a_yb_x\Delta\Delta a\, d\vec x- \int_{\mathbb T^2}a_{xy}\Delta b\Delta\Delta a\, d\vec x
- \int_{\mathbb T^2}a_{y}\Delta b\Delta\Delta a_x\, d\vec x\\
&+2 \int_{\mathbb T^2}a_{xy}b_{xx}\Delta\Delta a\, d\vec x+2 \int_{\mathbb T^2}a_{yy}b_{xy}\Delta\Delta a\, d\vec x\\
=:&\ M_{11}+M_{12}+M_{13}+M_{14}+M_{15},
\end{split}
\end{equation}
\begin{equation}\notag
\begin{split}
M_2=&-\int_{\mathbb T^2}\Delta a_xb_y\Delta\Delta a\, d\vec x- \int_{\mathbb T^2}a_x\Delta b_y\Delta\Delta a\, d\vec x\\
&-2 \int_{\mathbb T^2}a_{xx}b_{xy}\Delta\Delta a\, d\vec x-2 \int_{\mathbb T^2}a_{xy}b_{yy}\Delta\Delta a\, d\vec x\\
=&-\int_{\mathbb T^2}\Delta a_xb_y\Delta\Delta a\, d\vec x+ \int_{\mathbb T^2}a_{xy}\Delta b\Delta\Delta a\, d\vec x
+ \int_{\mathbb T^2}a_{x}\Delta b\Delta\Delta a_y\, d\vec x\\
&-2 \int_{\mathbb T^2}a_{xx}b_{xy}\Delta\Delta a\, d\vec x-2 \int_{\mathbb T^2}a_{xy}b_{yy}\Delta\Delta a\, d\vec x\\
=:&\ M_{21}+M_{22}+M_{23}+M_{24}+M_{25},
\end{split}
\end{equation}
and 
\begin{equation}\notag
\begin{split}
M_3=& \int_{\mathbb T^2}\Delta\left(a_y\Delta a_x\right)\Delta b\, d\vec x\\
=&\int_{\mathbb T^2}\Delta a_y \Delta a_x\Delta b\, d\vec x+ \int_{\mathbb T^2}a_y\Delta \Delta a_x\Delta b\, d\vec x\\
&+2 \int_{\mathbb T^2}a_{xy}\Delta a_{xx}\Delta b\, d\vec x+2 \int_{\mathbb T^2}a_{yy}\Delta a_{xy}\Delta b\, d\vec x\\
=:&\ M_{31}+M_{32}+M_{33}+M_{34},
\end{split}
\end{equation}
\begin{equation}\notag
\begin{split}
M_4=& -\int_{\mathbb T^2}\Delta\left(a_x\Delta a_y\right)\Delta b\, d\vec x\\
=&-\int_{\mathbb T^2}\Delta a_x \Delta a_y\Delta b\, d\vec x- \int_{\mathbb T^2}a_x\Delta \Delta a_y\Delta b\, d\vec x\\
&-2 \int_{\mathbb T^2}a_{xx}\Delta a_{xy}\Delta b\, d\vec x-2 \int_{\mathbb T^2}a_{xy}\Delta a_{yy}\Delta b\, d\vec x\\
=:&\ M_{41}+M_{42}+M_{43}+M_{44}.
\end{split}
\end{equation}
Again we observe the following cancellations,
\begin{equation}\notag
M_{12}+M_{22}=0, \ \ M_{13}+M_{32}=0, \ \ M_{23}+M_{42}=0, \ \ M_{31}+M_{41}=0.
\end{equation}
Hence we have
\begin{equation}\label{est-js1}
\begin{split}
&\ M_1+M_2+M_3+M_4\\
=&\  M_{11}+M_{14}+M_{15}+M_{21}+M_{24}+M_{25}+M_{33}+M_{34}+M_{43}+M_{44}.
\end{split}
\end{equation}
It follows from H\"older's, Gagliardo-Nirenberg's and Young's inequalities that
\begin{equation}\label{est-js2}
\begin{split}
|M_{11}|+|M_{21}|\leq &\ C\|\Delta\Delta a\|_{L^2} \|\nabla\Delta a\|_{L^{\frac{2r}{r-2}}} \|\nabla b\|_{L^r}\\
\leq &\ C\|\Delta\Delta a\|_{L^2}^{\frac{r+2}{r}} \|\nabla\Delta a\|_{L^2}^{\frac{r-2}{r}} \|\nabla b\|_{L^r}\\
\leq&\ \frac14\int_{\mathbb T^2}|\Delta \Delta a|^2\, d\vec x+C \|\nabla b\|_{L^r}^{\frac{2r}{r-2}}\|\nabla\Delta a\|_{L^2}^2
\end{split}
\end{equation}
By H\"older's inequality, Sobolev embedding and Young's inequality we obtain
\begin{equation}\notag
\begin{split}
|M_{14}|\leq &\ 2\|a_{xy}\|_{L^\infty(\mathbb T^2)}\|\Delta \Delta a\|_{L^2(\mathbb T^2)}\|b_{xx}\|_{L^2(\mathbb T^2)}\\
\leq &\ C\|\nabla \Delta a\|_{L^2(\mathbb T^2)}\|\Delta \Delta a\|_{L^2(\mathbb T^2)}\|b_{xx}\|_{L^2(\mathbb T^2)}\\
\leq &\ \frac{1}{32} \|\Delta \Delta a\|_{L^2(\mathbb T^2)}^2+C\|\nabla \Delta a\|_{L^2(\mathbb T^2)}^2\|b_{xx}\|_{L^2(\mathbb T^2)}^2
\end{split}
\end{equation}
for a constant $C$ that depends only on the size of the torus.
It is not hard to see that the other terms in (\ref{est-js1}) can be estimated in a similar way as that of $M_{14}$. Thus we have
\begin{equation}\label{est-js3}
\begin{split}
&\ |M_{14}|+|M_{15}|+|M_{24}|+|M_{25}|+|M_{33}|+|M_{34}|+|M_{43}|+|M_{44}|\\
\leq &\ \frac{1}{4} \|\Delta \Delta a\|_{L^2(\mathbb T^2)}^2+C\|\nabla \Delta a\|_{L^2(\mathbb T^2)}^2\|\nabla\nabla b\|_{L^2(\mathbb T^2)}^2.
\end{split}
\end{equation}

%{\color{blue} Question: we probably need to show that $b_{xy}\in L^2(0,T; L^2(\mathbb T^2))$ as well. }

Combing (\ref{est-ab8}) and (\ref{est-js1})-(\ref{est-js3}) we deduce 
\begin{equation}\label{est-ab9}
\begin{split}
&\frac{d}{dt} \int_{\mathbb T^2}|\nabla \Delta a|^2+|\Delta b|^2\, d\vec x+\mu \int_{\mathbb T^2}|\Delta \Delta a|^2+|\nabla\Delta b|^2\, d\vec x\\
\leq&\ C\left(\|\nabla b\|_{L^2}^{\frac{2r}{r-2}}+ \|\nabla\nabla b\|_{L^2(\mathbb T^2)}^2\right)\|\nabla \Delta a\|_{L^2(\mathbb T^2)}^2.
\end{split}
\end{equation}
Applying Gr\"onwall's inequality to (\ref{est-ab9}), combined with the assumption (\ref{criterion}) and the estimate (\ref{est-ab5}), we claim that 
\begin{equation}\label{est-ab10}
\begin{split}
& a\in L^\infty(0,T; H^3(\mathbb T^2))\cap L^2(0,T; H^4(\mathbb T^2)),\\
& b\in L^\infty(0,T; H^2(\mathbb T^2))\cap L^2(0,T; H^3(\mathbb T^2)).
\end{split}
\end{equation}

Such iteration process can be carried on to show that $(a(t), b(t))$ is smooth under the assumption (\ref{criterion}).

\begin{Remark}
We note that the diffusion term of $\Delta b$ is not used in the estimates above. 
%It indicates the diffusion in the $\vec e_z$ direction is not necessary in these results. 
\end{Remark}

\bigskip

\section*{Acknowledgement}
M. Dai is grateful for the hospitality of Leipzig University, Max Planck Institute at Leipzig, and Princeton University where part of the work was performed.

%\Endrefs
\end{document}